\documentclass[10pt,twoside]{amsart}
\usepackage[utf8]{inputenc}
\usepackage{default}
\usepackage{}
\usepackage{enumerate}
\usepackage{pgfpages}
\usepackage[leqno]{mathtools}
\usepackage{times}
\usepackage{amssymb}
\usepackage{amsmath}
\usepackage{textcomp}
 
\newcommand{\Pp}{\mbox{$\mathbb P$}}
\newcommand{\Cc}{\mbox{$\mathbb C$}}
\newcommand{\Zz}{\mbox{$\mathbb Z$}}
\DeclareMathOperator{\rank}{rank}

\newtheorem{thm}{Theorem}

\newtheorem{lem}{Lemma}[section]
\newtheorem{con}{Conjecture}[section]
\newtheorem{prop}[lem]{Proposition}

\newtheorem{rem}{Remark}

\usepackage{color}      

\title{ON A CONJECTURE ON LINEAR SYSTEMS}
\author{Sonica Anand}
\address{Indian Institute of Science Education and Research, Knowledge City, Sector-81, Mohali 140306 INDIA.}
\email{sonicaanand@iisermohali.ac.in}

\date{}
\subjclass[14H]{14H45,14H51}

\begin{document}

\maketitle 

\begin{abstract}
In a remark to Green's conjecture \cite{Gre84a}, Paranjape and Ramanan analyzed the vector bundle $E$ which is the pullback by the canonical map of the universal quotient bundle $T_{\Pp^{g-1}}(-1)$ on $\Pp^{g-1}$ \cite{PR88} and stated a more general conjecture \cite{HPR92} and proved it for the curves with Clifford Index $1$ (trigonal and plane quintics). In this paper, we state the conjecture for general linear systems and obtain results for the case of hyperelliptic curves.
\end{abstract}

\section{Introduction}\label{introduction}
Let $C$ be a smooth projective curve of genus $g \geq 2$ over a field $k$ and let $K$ be the canonical line bundle on $C$. In \cite{Gre84a}, Green made a conjecture which relates two aspects Koszul cohomology (an algebraic aspect) and Clifford Index $\gamma_C$ (a geometric aspect) of a curve. This conjecture \cite{Gre84a} is equivalent to the following \cite{PR88}: Let $E_K$ be the pullback by the canonical map $\Phi : C \rightarrow \Pp^{g-1}$ of universal quotient bundle on $\Pp^{g-1}$ . Then the map $\wedge ^i \Gamma(C, E_K) \rightarrow \Gamma(C, \wedge^i E_K)$ is surjective $\forall$ $i \leq \gamma_C$. Paranjape and Ramanan \cite{PR88} studied the vector bundle $E_K$ (stability properties). They also proved that all sections of $\wedge^i E_K$ which are locally decomposable are in the image of $\wedge^i \Gamma(E_K)  $ $\forall$ $i \leq \gamma_C$. Let $\sum_{i,K}$ be the cone of locally decomposable sections of $\wedge^iE_K$. In \cite{HPR92}, Hulek, Paranjape and Ramanan stated a conjecture.
\begin{con} \label{con1} $\sum_{i,K}$ spans $\Gamma(\wedge^i E_K)$ $\forall$ $i$ and for all curves.
\end{con}
This is stronger than Green's conjecture. They proved it for curves with Clifford index $1$ (trigonal curves and plane quintics). Conjecture \ref{con1} is trivial in case of hyperelliptic curves, since $E_K$ is the $(g-1)$- fold direct sum of the hyperelliptic line bundle. Vector bundle $E_K$ is semi-stable (even stable if $C$ is not hyperelliptic). In a remark to conjecture made in \cite{HPR92}, Eusen and Schreyer \cite{SE12} asked a more general question, whether $\Gamma(\wedge ^i N)$ is spanned by locally decomposable sections holds for every (stable) globally generated vector bundle $N$ on every curve $C$. They
gave counter examples to this more general question \cite{SE12}. By broadening our view point, in this paper we state a conjecture for general linear systems. Let $C$ be a smooth curve of genus $g \geq 2$ and let $L$ be a globally generated line bundle on $C$. The evaluation map gives rise to an exact sequence

\[ 0 \rightarrow E^* \rightarrow \Gamma(L)_C \rightarrow L \rightarrow 0\]
where $E^*$ is locally free of rank $h^0(L)-1$. Let $\sum_i$ be the cone of locally decomposable sections of $\wedge^iE$. We state
\begin{con}\label{con2} $\sum_i$ spans $\Gamma(\wedge^iE)$ $\forall$ $i$ and for all curves.\\
 \end{con}
In this paper we prove Conjecture \ref{con2} in case of hyperelliptic curves for the line bundles with degree large enough.\\

\begin{thm}\label{thm}
Let $C$ be a smooth hyperelliptic curve of genus $g \geq 2$ and let $L$ be a globally generated line bundle on $C$ of degree $d \geq 2g+1$ such that $H^1(L \otimes T^{-2})=0$, where $T$ is the hyperelliptic line bundle on $C$.  The evaluation map\index{Evaluation map} gives rise to an exact sequence 
   \begin{equation}\label{eval}
  0 \rightarrow E^* \rightarrow \Gamma(L)_C \rightarrow L \rightarrow 0
  \end{equation}
 where $E^*$ is locally free of $\rank h^0(L)-1$. Let $\sum_i$ be the cone of locally decomposable sections of $\wedge^iE$. Then
$\sum_i$ spans $\Gamma(\wedge^iE)$ $\forall$ $i$.
\end{thm}

This work has been done as a part of my Ph.D thesis under the guidance of Prof. Kapil H. Paranjape. I gratefully acknowledge IISER Mohali, the host institution and University Grant Commission, India for the financial support during this period.

\section{Geometry of the Hyperelliptic Curve}
 
%
 
 
Since $C$ is hyperelliptic of genus $g \geq 2$. Thus $g^1_2$ on $C$ is unique.
Let $\pi : C \rightarrow \Pp^1$ be the associated 2- sheeted covering, $T:= \pi^* \mathcal{O}_{\Pp^1}(1)$ is the unique $g^1_2$.\\
Consider the rank $2$ vector bundle $W$ on $\Pp^1$, where $W:= \pi_*L$. \\
Since,
\[ \chi(L)=\chi(\pi_*L)\]
So, 
\[d+1-g= rk(W) ( \frac{deg W}{rk W}+1-g_{\Pp^1}) = 2(\frac{deg W}{2}+1)\]
which gives \[deg W=d-g-1\]
Thus,
\begin{eqnarray*}
det W & \cong & \mathcal{O}_{\Pp^1}(d-g-1)\\
W & \cong & W^*(d-g-1)\\
\end{eqnarray*}

Since $\deg W = d-g-1 $, thus there is a unique integer $ x \leq \frac{d-g-1}{2}$ such that 
\begin{equation}\label{secW}
W \cong \mathcal{O}_{\Pp^1}(x) \bigoplus \mathcal{O}_{\Pp^1}(d-g-1-x)
\end{equation}

\begin{rem}
 \begin{enumerate}[(i)]
  \item $x$ is the least integer such that 
  \[ H^1(W(-2-x)) = H^1(\pi_*L(-2-x)) = H^1(L \otimes T^{-2-x}) \neq 0 \]
  In particular, this implies that 
  $\deg (L \otimes T^{-2-x}) \leq 2g-2 $ and thus we have 
  \[ \frac{d-2g-2}{2} \leq x \leq \frac{d-g-1}{2} \]
  \item Since $H^1(L \otimes T^{-2})=0$, thus $x > 0$, so we have
  \[ max.\{1,\frac{d-2g-2}{2} \} \leq x \leq \frac{d-g-1}{2} \]
  which implies both $W$ and $W(-1)$ is globally generated.
  \item Also, $H^1(L \otimes T^{-2})=0$ implies $H^1(L)=0$, thus by Riemann- Roch theorem, we have
  \begin{equation}\label{dim L}
 h^0(L)= d-g+1
 \end{equation} 
and $ \rank (E) = h^0(L)-1 = d-g \geq 3 $ ~~~~     (since $d \geq 2g+1$ and $g \geq 2$)
\item 
 We have 
\begin{eqnarray*}
\Gamma(W(-1)) & \cong & \Gamma(\pi_*(L \otimes T^{-1}))\\
 & \cong & \Gamma(L \otimes T^{-1})\\
\end{eqnarray*}
$H^1(L \otimes T^{-2})=0$ gives $H^1(L \otimes T^{-1})=0$\\
Thus, by Riemann-Roch theorem, we have \\
\[ h^0(L \otimes T^{-1})=d-g-1 \]
i.e., we have 
\begin{equation}\label{dimW(-1)}
h^0(W(-1))=d-g-1
\end{equation}        
\end{enumerate}
\end{rem}

\bigskip

Since $W(-1)$ is globally generated and $\Gamma(W) \cong \Gamma(L)$, so we have a surjection
\[ \Gamma(L)_{\Pp^1} \rightarrow W \rightarrow 0,\]
which is an isomorphism for sections. Since $W(-1)$ is globally generated bundle on $\Pp^1$, $W$ is very ample,\index{Very ample} i.e., we get an inclusion 
\begin{equation}\label{inclusion_projective}
 \Pp(W^*) \hookrightarrow \Pp(\Gamma(L)^*))=:\Pp.
\end{equation}
Also we have a surjection
\[ \pi^*W \rightarrow L \rightarrow 0\]
In other words, we have a subbundle of $\pi^*(W^*)$ that is isomorphic to $L^{-1}$. This gives a morphism from $C$ to $\Pp(W^*)$ with the property that the pullback of $\mathcal{O}_W(1)$ to $C$ is $L$. Also the composite of this morphism with the projection $p: \Pp(W^*) \rightarrow \Pp^1$ is $\pi$. Since the induced map
\[\Gamma(\Pp(W^*), \mathcal{O}_W(1)) \cong \Gamma(\Pp^1, W) \rightarrow \Gamma(C, L)\] is an isomorphism. Thus $C$ is actually embedded in $\Pp(W^*)$. Let us denote the image of $\Pp(W^*)$ in $\Pp$ by $S$.
We return to the ruled surface $p: \Pp(W^*) \rightarrow \Pp^1$. By (\ref{inclusion_projective}), there is an embedding $\Pp(W^*) \subset \Pp = \Pp(\Gamma(L)^*)$ with hyperplane section $\tau = \mathcal{O}_W(1)$.
Note that $\tau^2 = \deg W = d-g-1$. Since $C$ is a secant ($2$-section) of $\Pp(W^*)$, its class is of the form $ \mathcal{O}_W(2) \otimes p^* \mathcal{O}_{\Pp^1}(m)$. To compute $m$, we note that $d = C. \tau = 2 \tau^2 +m $.
Thus $m = 2g-d+2$.\\
Altogether, we have the following proposition\\
\begin{prop}\label{inclusion} There are inclusions
\[ C \subset S \subset \Pp\] with the following properties:\\
\begin{enumerate}[(i)]
 \item the restriction of $\mathcal{O}_{\Pp}(1)$ to $S$ is $\mathcal{O}_W(1)$;
 \item the restriction of $\mathcal{O}_{\Pp}(1)$ to $C$ is $L$;
 \item both restrictions induce isomorphisms of the corresponding linear systems;
 \item the divisor class on $S$ defined by $C$ is $\mathcal{O}_W(2) \otimes p^* \mathcal{O}_{\Pp^1}(-d+2g+2)$.
\end{enumerate}
\end{prop}

\textbf{Notation:} We will use the notation $U$ for the vector space $\Gamma(L \otimes T^{-1})$, i.e. we have
\begin{equation}\label{U}
\Gamma (W(-1)) \cong U
\end{equation}
and by (\ref{dimW(-1)}), we have
\begin{equation}\label{dimU}
 \dim U = d-g-1
\end{equation}


\section{Computation of dimensions}

In order to prove the conjecture, we want to relate the sections of $\wedge^iE$ to the sections of a suitable vector bundle on $\Pp^1$
\begin{lem}\label{EvalF} Let $F$ be a vector bundle on $\Pp^1$ that is globally generated. Then the evaluation sequence is 
\[ 0 \rightarrow \Gamma(F(-1))\otimes \mathcal{O}_{\Pp^1}(-1)\rightarrow \Gamma(F)_{\Pp^1}\rightarrow F \rightarrow 0 \]
\end{lem}
\textit{Proof}: $F$ is a sum of line bundles of degree $\geq 0$. Thus remains to check for line bundles, which is easy.\\
\qed

We want to apply this lemma to $W$ .\\
\begin{equation}\label{eval W}
0 \rightarrow \Gamma(W(-1))\otimes \mathcal{O}_{\Pp^1}(-1)\rightarrow \Gamma(W)_{\Pp^1}\rightarrow W \rightarrow 0
\end{equation}
Pulling back the evaluation sequence for $W$ on $\Pp^1$ to $C$ and using (\ref{U}) and the fact that $\Gamma(\pi_*L)\cong \Gamma(L)$, we get
\begin{equation}\label{pullW}
0\rightarrow U\otimes T^{-1}\rightarrow \Gamma(L)_C \rightarrow \pi^*W \rightarrow 0
\end{equation}
Also, we have a surjective map $\pi^*W \rightarrow L \rightarrow 0$, Let $Y$ be the kernel of $\pi^*W \rightarrow L \rightarrow 0$\\
i.e. we have \\
\begin{equation*}
0 \rightarrow Y \rightarrow \pi^*W \rightarrow L \rightarrow 0
\end{equation*}
\begin{eqnarray*}
\wedge ^2(\pi^*W) & \cong & Y \otimes L \\
\pi^*(\wedge^2 W) & \cong & Y \otimes L \\
\pi^*(\mathcal{O}_{\Pp^1}(d-g-1)) & \cong & Y \otimes L \\
 T^{d-g-1} & \cong & Y \otimes L \\
Y & \cong & L^{-1} \otimes T^{d-g-1}
\end{eqnarray*}
Thus, we have\\
\begin{equation}\label{ker}
0 \rightarrow L^{-1} \otimes T^{d-g-1} \rightarrow \pi^*W \rightarrow L \rightarrow 0
\end{equation}

we get a following commutative diagram\\

\begin{figure}
\centering
\includegraphics[scale= 0.8]{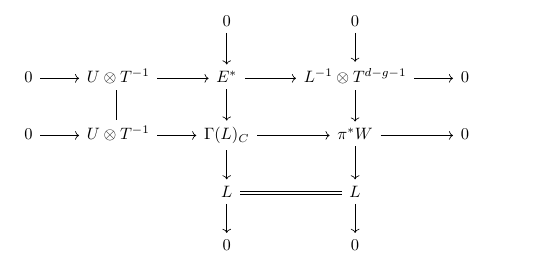}\\
\caption{}\label{fig1}
\end{figure}
\newpage
where the left vertical map is the evaluation map.

  Dualise the diagram in figure (\ref{fig1}), we get
  \begin{figure}[h]
\centering
\includegraphics[scale= 0.8]{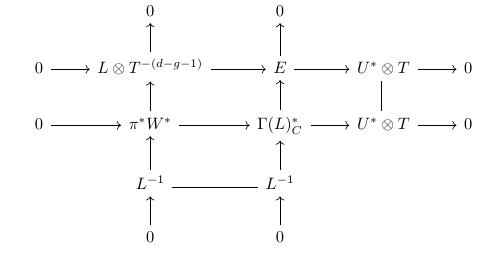}\\
\caption{}\label{fig2}
\end{figure}

The first line of commutative diagram in figure (\ref{fig2}) gives rise to an exact sequence
\begin{equation}\label{extE}
0 \rightarrow L \otimes T^{-(d-g-1)} \otimes \wedge^{i-1} U^*
\otimes T^{i-1} \rightarrow \wedge^i E \rightarrow \wedge^i U^*
\otimes T^{i} \rightarrow 0
\end{equation}

Since, $\pi^*W^*$ is a rank $2$ bundle we only get a filtration consisting of the following two exact sequences:

\begin{equation}\label{extL1}
0 \rightarrow \wedge^2 \pi^*W^* \otimes \wedge^{i-2} U^* \otimes T^{i-2} \rightarrow \wedge^i \Gamma(L)^*_C \rightarrow L_i \rightarrow 0
\end{equation}

\begin{equation}\label{extL2}
0 \rightarrow \pi^*W^* \otimes \wedge^{i-1}U^* \otimes T^{i-1} \rightarrow L_i \rightarrow \wedge^{i}U^* \otimes T^{i} \rightarrow 0
\end{equation}

Since the second horizontal sequence of diagram in figure (\ref{fig2}) is the pullback via $\pi$ of the dual of the sequence (\ref{eval W}). Thus both the above sequences come from $\Pp^1$, i.e. there exists a vector bundle $L'_i$ on $\Pp^1$, such that $L_i = \pi^*L'_i$ and the sequences
\begin{equation}\label{extLa}
0 \rightarrow \wedge^2 W^* \otimes \wedge^{i-2} U^* \otimes \mathcal{O}(i-2) \rightarrow \wedge^i \Gamma(L)^*_{\Pp^1} \rightarrow L'_i \rightarrow 0
\end{equation}

\begin{equation}\label{extLb}
0 \rightarrow W^* \otimes \wedge^{i-1}U^* \otimes \mathcal{O}(i-1) \rightarrow L'_i \rightarrow \wedge^{i}U^* \otimes \mathcal{O}(i) \rightarrow 0
\end{equation}

are such that (\ref{extL1}) and(\ref{extL2}) are pullback of (\ref{extLa}) and (\ref{extLb}) respectively.\\
Dualising (\ref{eval}), we have\\
\[ 0 \rightarrow L^{-1} \rightarrow \Gamma(L)_C^* \rightarrow E \rightarrow 0 \]
Thus the map $\wedge^i \Gamma(L)^*_C \rightarrow \wedge^i E$ is surjective.\\
Also we have $\wedge^i \Gamma(L)^*_C \rightarrow L_i \rightarrow 0$. The maps $\wedge^i \Gamma(L)^*_C \rightarrow \wedge^i E \rightarrow 0 $ factors through $L_i = \pi^*L'_i$. Thus we get the following commutative diagram with exact rows and columns:

\begin{figure}[h]
\centering
\includegraphics[scale= 0.8]{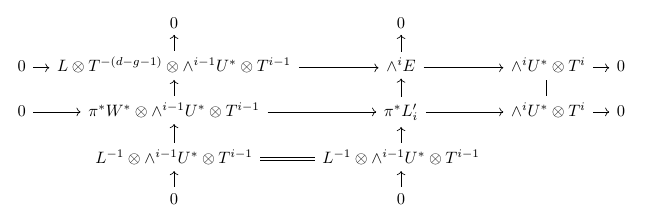}\\
\caption{}\label{fig3}
\end{figure}


where the top horizontal sequence is (\ref{extE}), middle horizontal sequence is (\ref{extL2}), left vertical sequence is obtained by dualising (\ref{ker}) and tensoring it with $\wedge^{i-1}U^* \otimes T^{i-1}$.

Let us compute the dimensions of the spaces $\Gamma(L'_i)$ and $\Gamma(\wedge^i E)$ for $ i \leq d-g$ ($= \rank E$)\\

\begin{lem}\label{lem2} When $d \geq 2g+1$, we have \[\dim \Gamma(L'_i) = \bigl (\begin{smallmatrix} d-g+1 \\
i \end{smallmatrix}\bigr) +\bigl (\begin{smallmatrix} d-g-1 \\
i-2 \end{smallmatrix}\bigr)(d-i-g)\] for $i \leq d-g$
\end{lem}

\textit{Proof}: 
Consider (\ref{extLa})\\
\[0 \rightarrow \wedge^2 W^* \otimes \wedge^{i-2} U^* \otimes \mathcal{O}(i-2) \rightarrow \wedge^i \Gamma(L)^*_{\Pp^1} \rightarrow L'_i \rightarrow 0 \]
Since \[det W \cong \mathcal{O}_{\Pp^1}(d-g-1).\] Thus \[det W^* \cong \mathcal{O}_{\Pp^1}(-(d-g-1))\]
Thus, we get
\[0 \rightarrow \wedge^{i-2} U^* \otimes \mathcal{O}(i-d+g-1) \rightarrow \wedge^i \Gamma(L)^* \rightarrow L'_i \rightarrow 0\]
 Since, \[ i \leq d-g,  h^0(\mathcal{O}(i-d+g-1))=0\]
Therefore,
\[h^0(L'_i) = h^0(\wedge^i \Gamma(L)^*) + h^1(\wedge^{i-2} U^* \otimes \mathcal{O}(i-d+g-1))\]
By (\ref{dimU}), we have
\[dim U = d-g-1\] and by (\ref{dim L})
 \[h^0(L) = d-g+1\] 
\[h^1(\mathcal{O}(i-d+g-1))= d-i-g\]
Thus,
\[\dim \Gamma(L'_i) = \bigl (\begin{smallmatrix} d-g+1 \\
i \end{smallmatrix}\bigr) +\bigl (\begin{smallmatrix} d-g-1 \\
i-2 \end{smallmatrix}\bigr)(d-i-g).\]
\qed

\subsection{Syzygies of the curve}
The syzygies of canonically embedded curves were computed by Schreyer \cite{S86}. Based on the parallel idea, we compute the syzygies of the curve $C$. For this, let
 \[ R = \bigoplus_{i=1}^{\infty} \Gamma(C, L^i)\] be the homogeneous coordinate ring\index{Homogeneous coordinate ring} of $C$ w.r.t $L$ and  
 \[ S= Sym \Gamma(C,L)= \bigoplus_{n \geq 0}\Gamma(\mathcal{O}_{\Pp^{d-g}}(n)).\]
 Let
\begin{equation}\label{resolution_F}
 0 \rightarrow F_t \rightarrow \cdots \rightarrow F_0 \rightarrow R \rightarrow 0,
\end{equation}
be a minimal free resolution \index{Minimal free resolution} of the graded $S$- module $R$. Then
 $F_i = \bigoplus_j S(-j)^{\beta_{ij}} = \bigoplus_j M_{ij} \otimes S(-j)$, where $M_{ij}$ is a $k$- vector space of $\dim \beta_{ij}$ and $S(-j)$ is the free $S$- module with one generator in degree $j$.\\
The resolution (\ref{resolution_F}) is equivalent to the free resolution of $\mathcal{O}_C$ as an $\mathcal{O}_{\Pp^{d-g}}$- module:
\[ 0 \rightarrow \bigoplus_j \mathcal{O}(-j)^{\beta_{d-g-1,j}} \rightarrow \cdots \rightarrow \bigoplus_j \mathcal{O}(-j)^{\beta_{0,j}} \rightarrow \mathcal{O}_C \rightarrow 0 \]
To find this resolution, one starts with the exact sequence 
\[ 0 \rightarrow \mathcal{O}_S(-2 \tau + (d-2g-2)f) \rightarrow \mathcal{O}_S \rightarrow \mathcal{O}_C \rightarrow 0 \]
(see Proposition (\ref{inclusion})). The idea is to first resolve the sheaves $\mathcal{O}_S$ and $\mathcal{O}_S(-2 \tau + (d-2g-2)f)$ resp. as $\mathcal{O}_{\Pp^{d-g}}$ modules and then form a mapping cone \index{Mapping cone}.
The result turns out to be a minimal resolution of $\mathcal{O}_C$.\\
Firstly, we will recall from \cite{Eis05} the description of the syzygies of these sheaves.\\
Let $\xi= \mathcal{O}(e_1)\bigoplus \mathcal{O}(e_2) \bigoplus \cdots \bigoplus \mathcal{O}(e_s)$ be a locally free sheaf of $\rank$ $s$ on $\Pp^1$, and let $p_\xi : \Pp(\xi) \rightarrow \Pp^1$ denote the corresponding $\Pp^{s-1}$ bundle. A rational normal scroll\index{Rational normal scroll} $X$ of type $S(e_1, e_2, \cdots, e_s)$ with $e_1 \geq e_2 \geq \cdots \geq e_s \geq 0$ and \[ f=e_1 + e_2 + \cdots + e_s \geq 2\] is the image of $\Pp(\xi)$ in $\Pp^r = \Pp(H^0(\Pp(\xi),\mathcal{O}_{\Pp(\xi)}(1)))$:
\[ j: \Pp(\xi) \rightarrow X \subset \Pp^r, r= f+s-1 \]
The Picard group\index{Picard group} of $\Pp(\xi)$ is generated by the hyperplane class $H=[j^*\mathcal{O}_{\Pp^r}(1)]$ and the ruling $R=[p_{\xi}^*\mathcal{O}_{\Pp^1}(1)]$:
\[ Pic \Pp(\xi) = \Zz H \bigoplus \Zz R, \] the intersection product is given by 
\[ H^s = f,   H^{s-1}. R = 1,    R^2 = 0. \]
We recall from \cite{Eis05}, the description of the syzygies of the sheaves 
\[ \mathcal{O}_X(aH+bR):= j_*\mathcal{O}_{\Pp(\xi)}(aH+bR),~~~ a,b \in \Zz \]
regarded as $\mathcal{O}_{\Pp^r}$- modules, at least in case $b \geq -1$.
\medskip

Let \[ \Phi: F \rightarrow G \]
be a map of locally free sheaves of $\rank$ $f'$ and $g'$, $f' \geq g'$, respectively on a smooth variety $V$.
We recall from \cite{BE75} the family of complexes $\zeta^b, b\geq -1$ of locally free sheaves on $V$, which resolve the $b^{th}$- symmetric power of coker $\Phi$ under suitable hypothesis on $\Phi$.\\
Define the $\j^{th}$ term in the complex $\zeta^b$ by \\
\[\zeta^b_j=\left\{\begin{array}{lll}
 \wedge^j F \otimes S_{b-j}G, &  for  & 0 \leq j \leq b\\
 \wedge^{j+g'-1}F \otimes D_{j-b-1}G^* \otimes \wedge^{g'} G^*, & for & j \geq b+1\\
 \end{array}\right.\]
and differential 
\[ \zeta^b_j \rightarrow \zeta^b_{j-1} \]
by the multiplication with $\Phi \in H^0(V, F^* \otimes G)$ for $j \neq b+1$ and $\wedge^{g'} \Phi \in H^0(V, \wedge^{g'}F^* \otimes \wedge^{g'}G)$ for $j=b+1$ in the appropriate term of the exterior $(\wedge F)$, symmetric (S.G) or divided power (D.G) algebra.

\begin{prop}\cite{Eis05}
$\zeta^b(a)$ for $b \geq -1$ is the minimal resolution of $\mathcal{O}_X(aH+bR)$ as an $\mathcal{O}_{\Pp^r}$- module, where $\zeta^b(a)= \zeta^b \otimes \mathcal{O}_{\Pp^r}(a)$ 
\end{prop}

\subsection{Minimal Resolution of $\mathcal{O}_C$}
We have \[ C \subset S \subset \Pp= \Pp(\Gamma(C,L)^*) \]
$C$ is contained in a $2$- dimensional rational normal scroll $S$ of type $S(e_1, e_2)$ and degree $f= e_1 + e_2 = d-g-1 \geq 2$. \\
$C$ is a divisor of class \[ C \sim 2H - (f-(g+1))R  ~~~~~ on ~~ S \]
The mapping cone \cite{S86}
\[ \zeta^{f-(g+1)}(-2) \rightarrow \zeta^{0} \] is the minimal resolution of $\mathcal{O}_C$ as an $\mathcal{O}_{\Pp^{d-g}}$- module.\\
We consider \[ \Phi: F \otimes \mathcal{O}_{\Pp^{d-g}}(-1) \rightarrow G \otimes \mathcal{O}_{\Pp^{d-g}} \]
be the map of locally free sheaves,
where $F$ is a vector space of dimension $f=d-g-1$ and $G$ is a vector space of dimension $2$\\

Firstly, we will compute 
\[ \zeta^{f-(g+1)}(-2) = \zeta^{d-2g-2} \otimes \mathcal{O}(-2) \]

Now,
\[ \zeta^{d-2g-2}_j = \left\{\begin{array}{ll}
                           \wedge^j (F \otimes \mathcal{O}(-1)) \otimes S_{d-2g-2-j}(G \otimes \mathcal{O}), & 0 \leq j \leq d-2g-2 \\
                           \wedge^{j+1} (F \otimes \mathcal{O}(-1)) \otimes D_{j-d+2g+2-1}(G \otimes \mathcal{O})^* \otimes \wedge^2 (G \otimes \mathcal{O})^*, & j \geq d-2g-1 \\
                          \end{array}\right. \]
 Since $j+1$ can be at most $d-g-1$. Thus, we have 
 \[ \zeta^{d-2g-2}_j = \left\{\begin{array}{ll}
                           \wedge^j (F \otimes \mathcal{O}(-1)) \otimes S_{d-2g-2-j}(G \otimes \mathcal{O}), & 0 \leq j \leq d-2g-2 \\
                           \wedge^{j+1} (F \otimes \mathcal{O}(-1)) \otimes D_{j-d+2g+2-1}(G \otimes \mathcal{O})^* \otimes \wedge^2 (G \otimes \mathcal{O})^*, & d-2g-1 \leq j \leq d-g-2 \\
                          \end{array}\right. \]
\[ \zeta^{d-2g-2}_j = \left\{\begin{array}{ll}
                           \wedge^j (F \otimes \mathcal{O}(-1)) \otimes S_{d-2g-2-j}(G \otimes \mathcal{O}), & 0 \leq j \leq d-2g-2 \\
                           \wedge^{j} (F \otimes \mathcal{O}(-1)) \otimes D_{j-d+2g}(G \otimes \mathcal{O})^* \otimes \wedge^2 (G \otimes \mathcal{O})^*, & d-2g-1 \leq j \leq d-g-1 \\
                          \end{array}\right. \]


  Similarly we can compute $\zeta^0_j$, 
 \[ \zeta^0_j = \left\{\begin{array}{ll}
                       \wedge^j(F \otimes \mathcal{O}(-1)) \otimes S_{0-j}(G \otimes \mathcal{O})   &  j=0 \\
\wedge^{j+1}(F \otimes \mathcal{O}(-1)) \otimes D_{j-1}(G \otimes \mathcal{O})^* \otimes \wedge^2(G \otimes \mathcal{O})^* & 1 \leq j \leq d-g-2 \\ 
                              \end{array}\right. \]  

The minimal free resolution\index{Minimal free resolution} of $\mathcal{O}_C$ is 
\begin{equation}\label{resolution_C}
0 \rightarrow \mathcal{O}(-(d-g-1))^g \rightarrow \cdots \rightarrow \wedge^{j+1} (F \otimes \mathcal{O}(-1)) \otimes D_{j-d+2g+2-1}(G \otimes \mathcal{O})^* \otimes \wedge^2 (G \otimes \mathcal{O})^* \rightarrow \cdots \rightarrow
\end{equation}
\[ \wedge^j (F \otimes \mathcal{O}(-1)) \otimes S_{d-2g-2-j}(G \otimes \mathcal{O}) \rightarrow \cdots \rightarrow \wedge^{j+1}(F \otimes \mathcal{O}(-1)) \otimes D_{j-1}(G \otimes \mathcal{O})^* \otimes \wedge^2(G \otimes \mathcal{O})^* \]
\[\rightarrow \cdots \rightarrow \wedge^j(F \otimes \mathcal{O}(-1)) \otimes S_{0-j}(G \otimes \mathcal{O}) \rightarrow \mathcal{O}_C \rightarrow 0 \]



We will use this resolution to compute the $\dim \Gamma(\wedge^i E)$. For this, consider (\ref{resolution_F}), the minimal free resolution of $R$ and recall the results from Chapter $1$ section $2$ of the Ph.D. thesis (1992) of Prof. Kapil H. Paranjape (University of Bombay, Bombay, India), we have
\[ M_{p,p+q} = coker (\wedge^{p+1}V \otimes \Gamma(C, K \otimes L^{q-1}) \rightarrow \Gamma(C, \wedge^p E^* \otimes L^q)) \]
where $M_{p,p+q} = (Tor_p^S(\Cc, R))_{p+q}$\\
and $\dim (Tor_p^S(\Cc, R))_{p+q} = \beta_{p,p+q}$.\\
Since $H^1(L)=0$, so we have
\[ M_{p,p+2} \approx H^1(\wedge^{p+1}E^*\otimes L)\]
\[ M_{p,p+2}^* \approx H^0(\wedge^{p+1}E \otimes L^{-1} \otimes K) \]

\begin{lem}\label{lem3} When $d \geq 2g+1$, we have \[\dim \Gamma(\wedge^iE) = \bigl (\begin{smallmatrix} d-g+1 \\
i \end{smallmatrix}\bigr) +\bigl (\begin{smallmatrix} d-g-1 \\
i-2 \end{smallmatrix}\bigr)(d-i-g)\]  for $i \leq d-g$
\end{lem}
\textit{Proof}: 
Consider (\ref{extE})\\
\[ 0 \rightarrow L \otimes T^{-(d-g-1)} \otimes \wedge^{i-1} U^* \otimes T^{i-1}
\rightarrow \wedge^i E \rightarrow \wedge^i U^* \otimes T^{i} \rightarrow 0\]
i.e.
\[0 \rightarrow L \otimes T^{-(d-g-i)} \otimes \wedge^{i-1} U^* \rightarrow \wedge^iE \rightarrow \wedge^i U^* \otimes T^{i} \rightarrow 0\]
Thus,

\begin{eqnarray*}
h^0(\wedge^i E) & = & \bigl [h^0(L \otimes T^{-(d-g-i)})-h^1(L \otimes T^{-(d-g-i)}) \bigr] \bigl (\begin{smallmatrix} d-g-1 \\
i-1 \end{smallmatrix}\bigr) + \bigl (\begin{smallmatrix} d-g-1 \\
i \end{smallmatrix}\bigr) h^0(T^i)\\
& &- \bigl (\begin{smallmatrix} d-g-1 \\
i \end{smallmatrix}\bigr) h^1(T^i) +h^1(\wedge^i E) \\
 & = &\bigl (\begin{smallmatrix} d-g-1 \\
i-1 \end{smallmatrix}\bigr)(g-d+2i+1) +\bigl (\begin{smallmatrix} d-g-1 \\
i \end{smallmatrix}\bigr)(i+1)- \bigl (\begin{smallmatrix} d-g-1 \\
i \end{smallmatrix}\bigr)(g-i) + h^1(\wedge^i E)  \\
\end{eqnarray*}
Now, \begin{eqnarray*}
      H^1(\wedge^iE) & = & H^1(\wedge^{d-g-i}E^* \otimes L)~~~~~~~(\rank E = d-g)\\
      & = & H^0(\wedge^{d-g-i}E \otimes L^{-1} \otimes K)^*
     \end{eqnarray*}
Thus, \[ h^1(\wedge^i E)= h^0(\wedge^{d-g-i}E \otimes L^{-1} \otimes K)\]
Since \[ M_{p,p+2}^* \approx H^0(\wedge^{p+1}E \otimes L^{-1} \otimes K) \]
Thus, we have \[ M_{d-g-i-1,d-g-i+1}^* = H^0(\wedge^{d-g-i}E \otimes L^{-1} \otimes K)\]
We have $d \geq 2g+1$ and $i \leq d-g$ which implies $i \leq g$. Thus, for $j= (d-g-1)-i \geq d-2g-1$, we have from (\ref{resolution_C}),
\[ \dim M_{d-g-i-1,d-g-i+1}^* = \dim \wedge^{d-g-1-i} (F \otimes \mathcal{O}(-1)) \otimes D_{g-1-i}(G \otimes \mathcal{O})^* \otimes \wedge^2 (G \otimes \mathcal{O})^* = \bigl (\begin{smallmatrix} d-g-1 \\
i \end{smallmatrix}\bigr)(g-i) \]
i.e. \[ h^1(\wedge^i E)=\bigl (\begin{smallmatrix} d-g-1 \\
i \end{smallmatrix}\bigr)(g-i) \]

Thus
\begin{eqnarray*}
\dim \Gamma(\wedge^iE)& = &\bigl (\begin{smallmatrix} d-g-1 \\
i-1 \end{smallmatrix}\bigr)(g-d+2i+1) +\bigl (\begin{smallmatrix} d-g-1 \\
i \end{smallmatrix}\bigr)(i+1)\\
& = &\bigl (\begin{smallmatrix} d-g-1 \\
i-1 \end{smallmatrix}\bigr)(g-d+2i+1) +\bigl (\begin{smallmatrix} d-g-1 \\
i \end{smallmatrix}\bigr)(i+1)+ \bigl (\begin{smallmatrix} d-g+1 \\
i \end{smallmatrix}\bigr)- \bigl (\begin{smallmatrix} d-g+1 \\
i \end{smallmatrix}\bigr)\\
& = & \bigl (\begin{smallmatrix} d-g+1 \\
i \end{smallmatrix}\bigr) +\bigl (\begin{smallmatrix} d-g-1 \\
i-2 \end{smallmatrix}\bigr)(d-i-g)\\
\end{eqnarray*}
\qed


\begin{prop} \label{prop1} For $i \leq d-g$, the map $\Gamma(L_i') \rightarrow \Gamma(\wedge^i E)$ induced by diagram in figure (\ref{fig3}) is an isomorphism.\\
\end{prop}
\textit{Proof}: We get an exact commutative diagram\\

\medskip
\begin{figure}[h]
\centering
\includegraphics[scale= 0.8]{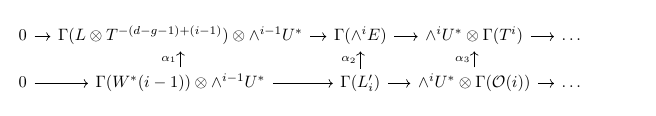}\\
\end{figure}
  
   \medskip

where $\alpha_i's$ are induced by diagram in figure (\ref{fig3}). 
\begin{itemize}
\item $\alpha_1$ is injective. \\
Consider the left vertical sequence of diagram in figure (\ref{fig3})\\
\begin{equation*}
0\rightarrow L^{-1} \otimes \wedge^{i-1} U^* \otimes T^{i-1} \rightarrow \pi^*W^* \otimes \wedge^{i-1} U^* \otimes T^{i-1} 
\end{equation*}
\begin{equation*}
\rightarrow L \otimes T^{-(d-g-1)} \otimes \wedge^{i-1} U^*
\otimes T^{i-1} \rightarrow 0
\end{equation*}\\

This gives rise to\\
\begin{equation*}
0\rightarrow \Gamma(L^{-1}\otimes T^{i-1}) \otimes \wedge^{i-1} U^*  \rightarrow \Gamma(W^*(i-1)) \otimes \wedge^{i-1} U^* 
\end{equation*}
\begin{equation*}
\underrightarrow{\alpha_1} \Gamma(L \otimes T^{-(d-g-1)+(i-1)}) \otimes \wedge^{i-1} U^*
 \rightarrow \ldots
\end{equation*}\\
Since
\begin{eqnarray*}
\Gamma(W^* \otimes \mathcal{O}_{\Pp^1}(i-1)) & \cong & \Gamma(W \otimes \mathcal{O}_{\Pp^1}(-d+g+i))\\
& \cong & \Gamma(\pi_*L \otimes \mathcal{O}_{\Pp^1}(-d+g+i))\\
&\cong & \Gamma(\pi_*(L \otimes \pi^*\mathcal{O}_{\Pp^1}(-d+g+i)))\\
& \cong & \Gamma(L \otimes \pi^*\mathcal{O}_{\Pp^1}(-d+g+i))\\
& \cong & \Gamma(L \otimes T^{-d+g+i})
\end{eqnarray*}

Thus $\alpha_1$ is injective.

\item $\alpha_3$ is injective by definition.\\
\end{itemize}
Hence $\alpha_2$ is injective and since both the spaces have the same dimension, the map $\Gamma(L_i') \rightarrow \Gamma(\wedge^i E)$ is an isomorphism.\\
\qed

\section{Construction of a subbundle of $E$}

We want to prove the conjecture for hyperelliptic curves of genus g. We shall first do this for $i=2$. The main point is to construct sufficiently many locally decomposable sections that are not globally decomposable.\\
Consider $p:\Pp(W^*) \rightarrow \Pp^1 $, the natural projection. For every $a \in \Pp^1$, the fibre $l_a = p^{-1}(a)$ is a secant of the curve $C$. Let $W^*_a$ be the fibre of $W^*$ at $a$. Since $W$ is globally generated thus we have $\Gamma(W)_{\Pp^1} \rightarrow W \rightarrow 0$. $W= \pi_*L$ and $\Gamma(\pi_*L) \cong \Gamma(L)$ thus we have $\Gamma(L)_{\Pp^1} \rightarrow W \rightarrow 0$, which gives $W^* \rightarrow \Gamma(L)^*_{\Pp^1}$ and we can identify $\Gamma(L)^*$ with $\Gamma(E)$. Also, we have 
\[0 \rightarrow E^* \rightarrow \Gamma(L)_C \rightarrow L \rightarrow 0,\] which gives $\Gamma(L)^*_C \rightarrow E$ i.e. a map $\Gamma(E)_C \rightarrow E$. Thus we get a map 
\[(W_a)^*_C \rightarrow \Gamma(E)_C \rightarrow E\] which is composite of the inclusion of $W_a^*$ in $\Gamma(E)$ and the evaluation map.

Let $F(a)$ be the subbundle of E generated by the image of $W_a^*$. A section of $\Gamma(E)$ is non-zero at every point of $C$ if it corresponds to a point of $\Pp(\Gamma(L)^*)$ not on the curve $C$, while a section corresponding to a point say $x \in C$ vanishes exactly at $x$.\\
Hence the map $W_a^* \rightarrow F(a)$ is an isomorphism outside $C \bigcap l_a$ but has rank $1$ over $C \bigcap l_a$.The induced map $\wedge^2 W_a^* \rightarrow \wedge^2 F(a)$ has simple zeros exactly over $C \bigcap l_a$. Hence $F(a)$ has rank $2$ and $\wedge^2 F(a)= T$.\\
The vector bundle $F(a)$ has $ W_a^*$ as its space of sections i.e. $\dim \Gamma(F(a))=2$. On the other hand $\dim \Gamma(\wedge^2 F(a))= \dim \Gamma(T)= 2$. Thus we get a $2$- dimensional subspace of $\Gamma(\wedge^2E)$ consisting of locally decomposable sections of which only the $1$- dimensional subspace $\wedge^2 \Gamma(F(a)) \subset \Gamma(\wedge^2F(a))$ consists of globally decomposable sections.\\
The next step is to globalise this construction, i.e. to vary the point $a$. We consider the graph inclusion $\Gamma \subset C \times \Pp^1$ given by the map $\pi$. This divisor belongs to the line bundle $p_1^*T \otimes p_2^* \mathcal{O} _{\Pp^1}(1)$, where $p_1$ and $p_2$ are the natural projections to $C$, resp. $\Pp^1$. the direct image by $p_2$ of the bundle morphism $ p_2^*W^* \rightarrow \Gamma(E)_{C \times \Pp^1}$ yields the map $W^* \rightarrow \Gamma(E)_{\Pp^1}$, and hence a map $\wedge ^2W^* \rightarrow \wedge^2 \Gamma(E)_{\Pp^1}$.\\
On the other hand the bundle homomorphism $p_2^* W^* \rightarrow \Gamma(E)_{C \times \Pp^1} \rightarrow p_1^*E$ fails to be injective precisely over $\Gamma$. Thus, we get a morphism\\
$p_2^*(\wedge^2W^*) \otimes \mathcal{O}(\Gamma) \rightarrow p_1^*(\wedge^2 E)$.
Taking direct image by $p_2$ gives a morphism\\
$\wedge^2 W^* \otimes \Gamma(T) \otimes \mathcal{O}(1) \rightarrow \Gamma(\wedge^2 E)_{\Pp^1}$. For every $a \in \Pp^1$ this induces a map $\Gamma(T) \rightarrow \Gamma(\wedge^2 E)_{\Pp^1}$ and this gives exactly the space of locally decomposable sections described above.\\

Altogether, we get a commutative diagram\\
\begin{figure}[h]
\centering
\includegraphics[scale= 0.8]{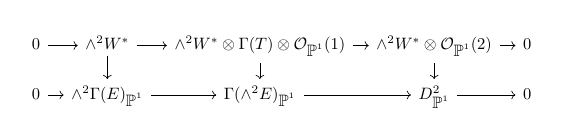}\\
\end{figure}


where $D^2_{\Pp^1} := \frac{\Gamma(\wedge^2E)}{\wedge^2\Gamma(E)}$\\
where the top horizontal row is the evaluation sequence for $\mathcal{O}_{\Pp^1}(1)$, which is 
\[0 \rightarrow \mathcal{O}_{\Pp^1}(-1) \rightarrow \Gamma(\mathcal{O}_{\Pp^1}(1))\otimes \mathcal{O}_{\Pp^1} \rightarrow \mathcal{O}_{\Pp^1}(1) \rightarrow 0\]
Tensoring with $\wedge^2W^* \otimes \mathcal{O}_{\Pp^1}(1)$, we get
\[ 0 \rightarrow \wedge^2W^* \rightarrow \wedge^2W^* \otimes \Gamma(T) \otimes \mathcal{O}_{\Pp^1}(1)  \rightarrow \wedge^2W^* \otimes \mathcal{O}_{\Pp^1}(2)  \rightarrow 0 \] and use $\Gamma(T) \cong \Gamma(\mathcal{O}_{\Pp^1})$.\\
We have to show that the locally decomposable sections constructed above together with $\wedge^2 \Gamma(E)$ generate $\Gamma(\wedge^2E)$. For this, we consider the map $\wedge^2W^* \otimes \mathcal{O}_{\Pp^1}(2) \rightarrow D^2_{\Pp^1}$. We want to show that this map is injective (as a bundle map) and that the resulting rational curve in $\Pp(D^2)$ is the rational normal curve of degree $d-g-3$ (recall that $\dim D^2 = d-g-2$). This is sufficient since the rational normal curve of degree n in $\Pp^n$ spans $\Pp^n$.

Our aim is to do this by entirely reducing the problem to computations on $\Pp^1$, resp. $\Pp^1 \times \Pp^1$.\\
\begin{lem}\cite{HPR92} \label{lem} Let $ \mathcal{O}_{\Pp^1}(-n) \rightarrow \Gamma(\mathcal{O}_{\Pp^1}(n))^*_{\Pp^1}$ be a non-zero $Sl_2(\Cc)$- equivariant morphism. Then this morphism defines an embedding of $\Pp^1$ into $\Pp(\Gamma(\mathcal{O}_{\Pp^1}(n)^*))$ as a rational normal curve of degree $n$.
\end{lem}

We return to the bundle $W$. Sequence (\ref{extLa}) gives for $i=2$ the following sequence:\\
\begin{equation}\label{Wi=2}
0 \rightarrow \wedge^2W^* \rightarrow \wedge^2\Gamma(W)^*_{\Pp^1} \rightarrow L_2' \rightarrow 0
\end{equation}

Consider $\Pp^1 \times \Pp^1$ together with projections $q_1$ and $q_2$ resp.\\
Taking pullback of (\ref{Wi=2}) via $q_1$ and $q_2$ resp., we get a map $q_2^*\wedge^2W^* \rightarrow q_1^*L_2'$ that vanishes along the diagonal $\bigtriangleup \subset \Pp^1 \times \Pp^1$.\\
Hence, we get a morphism 
\[ q_2^*\wedge^2W^* \otimes \mathcal{O}(\bigtriangleup)\rightarrow q_1^*L_2'\]
Applying $ q_{2*} $, we get a map
\[ \wedge^2W^* \otimes \Gamma(\mathcal{O}(1)) \otimes \mathcal{O}(1) \rightarrow \Gamma(L_2')_{\Pp^1}\]
This gives rise to a commutative diagram\\

\begin{figure}[h]
\centering
\includegraphics[scale= 0.8]{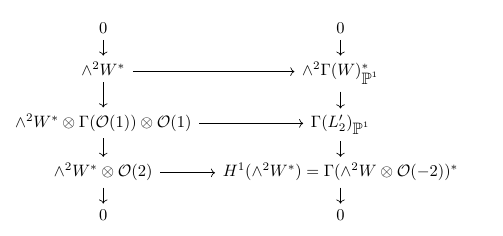}\\
\caption{}\label{fig6}
\end{figure}


where the left hand column is the Euler sequence on $\Pp^1$ twisted by $\wedge^2W^* \otimes \mathcal{O}(1)$, the right hand column comes from (\ref{Wi=2}) and the map $\wedge^2W^* \rightarrow \wedge^2 \Gamma(W)^*$ is the natural one. This diagram is $Sl_2(\Cc)$ equivariant, where $Sl_2(\Cc)$ acts on $\Pp^1$ in the usual way and on $\Pp^1 \times \Pp^1$ by the diagonal action. In particular the morphism $\wedge^2 W^* \otimes \mathcal{O}(2) \rightarrow \Gamma(\wedge^2 W \otimes \mathcal{O}(-2))^*$ is $Sl_2(\Cc)$ equivariant, by Lemma(\ref{lem}), it defines an embedding of $\Pp^1$ into $\Pp(\Gamma(\wedge^2 W \otimes \mathcal{O}(-2))^*)$ as a rational normal curve of degree $d-g-3$ .\\

\begin{lem}\cite{HPR92} \label{lem} Diagram in figure (\ref{fig6}) gives rise to a commutative and exact diagram
\end{lem}

\begin{figure}[h]
\centering
\includegraphics[scale= 0.75]{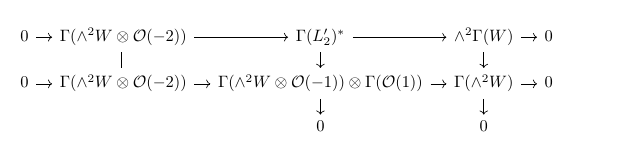}\\
\end{figure}

\newpage
 
 \begin{prop}\label{prop} $\Gamma(\wedge^2 E)$ is generated by locally decomposable sections.
 \end{prop}
 \textit{Proof}:
 We have constructed maps $p_2^*(\wedge^2W^*) \otimes \mathcal{O}(\Gamma) \rightarrow p_1^*(\wedge^2 E)$ on $ C \times \Pp^1$
 and $q_2^*\wedge^2W^* \otimes \mathcal{O}(\bigtriangleup)\rightarrow q_1^*L_2'$ on $\Pp^1 \times \Pp^1$. Consider the diagram
 
 \begin{figure}[h]
\centering
\includegraphics[scale= 0.8]{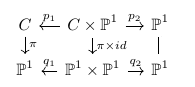}\\
\end{figure}
Pulling the morphism $q_2^*\wedge^2W^* \otimes \mathcal{O}(\bigtriangleup)\rightarrow q_1^*L_2'$ on $\Pp^1 \times \Pp^1$ back to $ C \times \Pp^1$, we get a morphism $p_2^*(\wedge^2W^*) \otimes \mathcal{O}(\Gamma) \rightarrow p_1^*(\pi^* L_2')$. By construction the diagram \\
\begin{figure}[h]
\centering
\includegraphics[scale= 0.8]{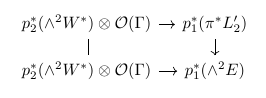}\\
\end{figure}
 commutes where the map  $p_1^*(\pi^* L_2') \rightarrow p_1^* \wedge^2E $ is the pullback via $p_1$ of the corresponding map in diagram in figure (\ref{fig3}).Pushing this down via  $\pi \times id$ to $\Pp^1 \times \Pp^1$ leads to the commutative diagram
 \begin{figure}[h]
\centering
\includegraphics[scale= 0.8]{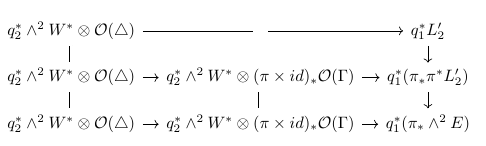}\\
\end{figure}
 
 Now taking $q_{2*}$ of the outermost square we get 
 \begin{figure}[h]
\centering
\includegraphics[scale= 0.8]{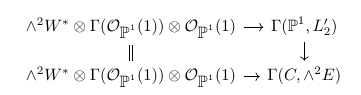}\\
\end{figure}
  
 where the right hand vertical map is an isomorphism from Proposition \ref{prop}.\\
 Thus in order to compute the diagram \\
 \begin{figure}[h]
\centering
\includegraphics[scale= 0.8]{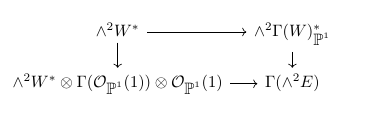}\\
\end{figure}

\newpage

  we can compute
  \begin{figure}[h]
\centering
\includegraphics[scale= 0.8]{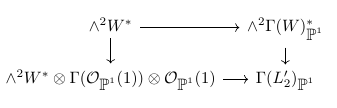}\\
\end{figure}
  
  and the result follows from Lemma \ref{lem}.
  \qed

  \section{proof of main result}
  Here we want to prove Theorem \ref{thm}.\\
%
We shall first show that for $2 \leq i \leq d-g$, there is a natural epimorphism
 \[\wedge^{i-2} \Gamma(W)^* \otimes \Gamma(L_2')_{\Pp^1} \rightarrow \Gamma(L_i')_{\Pp^1} \rightarrow 0\]
 Setting $i=2$ in (\ref{extLa}), we get
 \[0 \rightarrow \wedge^2W^* \rightarrow \wedge^2\Gamma(L)^*_{\Pp^1} \rightarrow L_2' \rightarrow 0\]
 Twisting with $\wedge^{i-2} \Gamma(W)^*$, we get an exact sequence \[0 \rightarrow \wedge^{i-2} \Gamma(W)^* \otimes \wedge^2W^* \rightarrow \wedge^{i-2} \Gamma(W)^* \otimes \wedge^2\Gamma(W)^* \rightarrow \wedge^{i-2} \Gamma(W)^* \otimes  L_2' \rightarrow 0\] (since $\Gamma(L)\cong \Gamma(W)$)
  Combining this with (\ref{extLa}), we get a diagram 
  \begin{figure}[h]
\centering
\includegraphics[scale= 0.85]{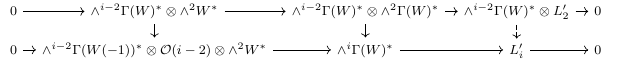}\\
\caption{}\label{fig14}
\end{figure}

 
 Here the middle vertical map is the canonical one and the left hand vertical map is given by taking $\wedge^{i-2}$ of the dual evaluation sequence 
 \[0 \rightarrow W^* \rightarrow \Gamma(W)^* \otimes \mathcal{O}_{\Pp^1} \rightarrow \Gamma(W(-1))^* \otimes \mathcal{O}_{\Pp^1}(1) \rightarrow 0 \]
 Taking $\wedge^{i-2}$ of the above sequence, we get
 \[0 \rightarrow \wedge^2W^* \otimes \wedge^{i-2} \Gamma(W(-1))^* \otimes \mathcal{O}_{\Pp^1}(i-2) \rightarrow \wedge^{i-2} \Gamma(W)^* \rightarrow F_{i-2} \rightarrow 0 \]
 \[0 \rightarrow F_{i-2} \rightarrow \wedge^{i-2} \Gamma(W)^* \otimes \mathcal{O}_{\Pp^1} \rightarrow \wedge^{i-2}\Gamma(W(-1))^* \otimes \mathcal{O}_{\Pp^1}(i-2) \rightarrow 0\]
 
 Tensoring above sequence with $\wedge ^2 W^*$, we get
 \[0 \rightarrow F_{i-2} \otimes \wedge ^2 W^* \rightarrow \wedge^{i-2} \Gamma(W)^* \otimes \mathcal{O}_{\Pp^1} \otimes \wedge ^2 W^* \rightarrow \wedge^{i-2}\Gamma(W(-1))^* \otimes \mathcal{O}_{\Pp^1}(i-2)\otimes \wedge ^2 W^* \rightarrow 0\]
 
 Taking the associated cohomology sequence of commutative diagram in figure (\ref{fig14}), we get the following commutative diagram:\\
 \begin{figure}[h]
\centering
\includegraphics[scale= 0.8]{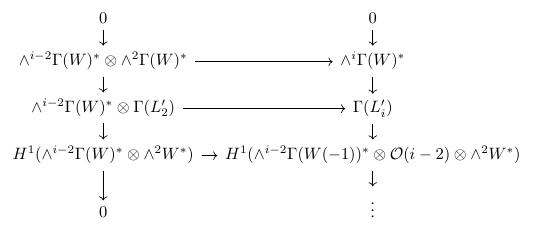}\\
\end{figure}

Here $W$ is a rank $2$ vector bundle on $\Pp^1$ of degree $d-g-1$.
\begin{eqnarray*}
det W & \cong & \mathcal{O}(d-g-1)\\
det W^* & \cong & \mathcal{O}(-d+g+1)\\
i.e. \wedge^2 W^* & \cong & \mathcal{O}(-d+g+1)
\end{eqnarray*}
Since $2 \leq d-g-1$ \\
Thus $\Gamma(\wedge^2 W) =0$\\

The top horizontal map is clearly surjective. The bottom horizontal map is surjective since $H^2$  vanishes on $\Pp^1$. By standard diagram chasing the middle horizontal map must be surjective thus giving our first claim.\\
By construction the natural diagram\\
\begin{figure}[h]
\centering
\includegraphics[scale= 0.8]{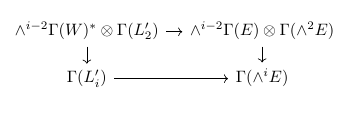}\\
\end{figure}

Since $\Gamma(W)^* \cong \Gamma(L)^*$ and $\Gamma(L)^*$ can be identified with $\Gamma(E)$. By proposition (\ref{prop}), the horizontal maps are isomorphisms. Hence the natural map $\wedge^{i-2} \Gamma(E) \otimes \Gamma(\wedge^2E) \rightarrow \Gamma(\wedge^i E)$ is surjective and our claim follows from Proposition (\ref{prop}).

\bibliographystyle{amsalpha}
\bibliography{Ref}

\providecommand{\bysame}{\leavevmode\hbox to3em{\hrulefill}\thinspace}
\providecommand{\MR}{\relax\ifhmode\unskip\space\fi MR }
\providecommand{\MRhref}[2]{%
  \href{http://www.ams.org/mathscinet-getitem?mr=#1}{#2}
}
\providecommand{\href}[2]{#2}
\begin{thebibliography}{HPR92}

\bibitem[BE75]{BE75}
David~A. Buchsbaum and David Eisenbud, \emph{Generic free resolutions and a
  family of generically perfect ideals}, Advances in Math. \textbf{18} (1975),
  no.~3, 245--301. \MR{0396528}

\bibitem[Eis05]{Eis05}
David Eisenbud, \emph{The geometry of syzygies}, Graduate Texts in Mathematics,
  vol. 229, Springer-Verlag, New York, 2005, A second course in commutative
  algebra and algebraic geometry. \MR{2103875}

\bibitem[ES12]{SE12}
Friedrich Eusen and Frank-Olaf Schreyer, \emph{A remark on a conjecture of
  {P}aranjape and {R}amanan}, Geometry and arithmetic, EMS Ser. Congr. Rep.,
  Eur. Math. Soc., Z\"urich, 2012, pp.~113--123. \MR{2987656}

\bibitem[Gre84]{Gre84a}
Mark~L. Green, \emph{Koszul cohomology and the geometry of projective
  varieties}, J. Differential Geom. \textbf{19} (1984), no.~1, 125--171.
  \MR{739785 (85e:14022)}

\bibitem[HPR92]{HPR92}
K.~Hulek, K.~Paranjape, and S.~Ramanan, \emph{On a conjecture on canonical
  curves}, J. Algebraic Geom. \textbf{1} (1992), no.~3, 335--359. \MR{1158621
  (93c:14029)}

\bibitem[PR88]{PR88}
Kapil Paranjape and S.~Ramanan, \emph{On the canonical ring of a curve},
  Algebraic geometry and commutative algebra, {V}ol.\ {II}, Kinokuniya, Tokyo,
  1988, pp.~503--516. \MR{977775 (90b:14024)}

\bibitem[Sch86]{S86}
Frank-Olaf Schreyer, \emph{Syzygies of canonical curves and special linear
  series}, Math. Ann. \textbf{275} (1986), no.~1, 105--137. \MR{849058}

\end{thebibliography}

\end{document}